\documentclass[10pt,oneside]{amsart} 

\usepackage{amssymb} 
\usepackage[mathscr]{eucal} 
\usepackage[matrix,arrow,tips,curve,ps]{xy}
\usepackage{amsmath} 
\usepackage{epsfig}
\usepackage{amscd}
\newtheorem{TEO}{Theorem}[section]

\theoremstyle{definition}

\def\OO{{\mathcal O}}

\newcommand\dual{\mathrel{\raise3pt\hbox{$\underline{\mathrm{\thinspace d
\thinspace}}$}}}

\newcommand\proj{\mathbb P}

\newcommand\Co{\mathbb C}

\def\Co{{\mathbb C}}

\newcommand{\ud}{\mathrm{d}} 
\newcommand{\LB}{\omega_C \otimes {A}}

\begin{document}

\footnote{
Partially supported by  PRIN  2009: "Moduli, strutture geometriche e loro applicazioni" and by INdAM (GNSAGA).
 }

\title{On the first Gaussian map for Prym-canonical line bundles}

\author[C. Barchielli]{Caterina Barchielli}
\address{Dipartimento di Matematica,
Universit\`a di Pavia,
via Ferrata 1, I-27100 Pavia, Italy } \email{{\tt
caterina.barchielli@tiscali.it}}

\author[P. Frediani]{Paola Frediani}
\address{ Dipartimento di Matematica, Universit\`a di Pavia,
via Ferrata 1, I-27100 Pavia, Italy } \email{{\tt
paola.frediani@unipv.it}}

\maketitle

\setlength{\parskip}{.1 in}

\begin{abstract}
We prove by degeneration to Prym-canonical binary curves that the first Gaussian map $\mu_A$ of the Prym canonical line bundle  $\omega_C \otimes A$ is surjective for the general point $[C,A] \in {\mathcal R}_g$ if $g \geq 12$, while it is injective if $g \leq 11$.
\end{abstract}

\section{Introduction}
Consider the moduli space ${\mathcal R}_g$ parametrizing isomorphism classes of pairs $[C,A]$, where $C$ is a smooth projective curve of genus $g$ and $A$ is a non trivial two torsion line bundle on $C$. 
In this paper we prove that for a general point $[C,A]$ in ${\mathcal R}_g$, the first Gaussian map of the Prym-canonical line bundle $\omega_C \otimes A$ is of maximal rank, i.e. it is surjective for $g \geq 12$ and injective for $g \leq 11$. 

We do it by degeneration to certain stable binary curves of genus $g$ embedded in ${\proj}^{g-2}$ by a Prym-canonical linear system $|\omega_C \otimes A|$ that have been constructed in \cite{cf} to prove the surjectivity of the second Gaussian map of such line bundles for the general element in ${\mathcal R}_g$ for $g \geq 20$.  
A stable binary curve of genus $g$ is the union of two smooth rational curves meeting transversally at $g+1$ distinct points and in \cite{ccm} Calabri, Ciliberto and Miranda have first used these  curves to prove by degeneration the maximality of the rank of the second Gaussian map of the canonical line bundle for a general curve. 

The first Gaussian (or Wahl) map of the canonical line bundle has been extensively studied by many authors. Here we only recall that  Wahl's theorem (\cite{wahl}, see also \cite{bm}) says that if  $C$ is a curve sitting on a K3 surface, then its first Wahl map is not surjective, while Ciliberto, Harris and Miranda have proved in \cite{chm} that the first Wahl map is surjective for the general curve of genus $g \geq 10$, $g \neq 11$ (see also \cite{voi}).

We also point out that in \cite{cv} the Gaussian map $\gamma^1_{\omega_C, \omega_C \otimes A}$ defined in section 3.1 (where $A$ is a 2-torsion line bundle) has been studied. The main result of \cite{cv} is the proof of the surjectivity of the map $\gamma^1_{\omega_C, \omega_C \otimes A}$ for the general point in ${\mathcal R}_g$ for $g \geq 12$, $g \neq 13, 19$.  Their method is to use degeneration to graph curves as in \cite{chm}. For dimension reasons one would expect the surjectivity of the map $\gamma^1_{\omega_C, \omega_C \otimes A}$ for the general point in ${\mathcal R}_g$, as long as $g \geq 8$. The restriction $g \geq 12$, $g \neq 13$ that they impose is only due to the combinatorics of the graph curves. They also observe that a similar result to the above mentioned theorem of Wahl on the non surjectivity of the first Wahl map for hyperplane sections of K3 surfaces holds, namely, if $C$ is a hyperplane section of an Enriques surface, then the map $\gamma^1_{\omega_C, \omega_C \otimes A}$ is not surjective, where $A$ is the line bundle of order 2 induced from the canonical bundle of the surface.

Here we focus instead on the study of the first Gaussian map $\mu_A: \Lambda^2 H^0(C, \omega_C \otimes A) \to H^0(C, \omega_C^{\otimes 3})$, where $\mu_A$ is the Gaussian map $\gamma^1_{ \omega_C \otimes A,  \omega_C \otimes A}$ defined in section 3.1. 

The main result of this paper says that the first Gaussian map $\mu_A$ of a general Prym-canonical binary curve of genus $g$ is surjective for $g \geq 12$. In the case $g = 12$ we verify the surjectivity of $\mu_A$ on explicit Prym-canonical binary curves using Maple and similar computations with Maple also give the injectivity of $\mu_A$ for a general Prym-canonical binary curve of  genus $g \leq 11$. Then, to prove surjectivity for higher genus $g$ we proceed by induction on $g$ in a similar way as it has been done in \cite{ccm} and in \cite{cf}.
Maximality of the rank for the general element in ${\mathcal R}_g$ follows  then by semicontinuity. 

This implies that the locus of points $[C,A]$ in ${\mathcal R}_g$ for which the map $\mu_A$ is not surjective is a proper subscheme of ${\mathcal R}_g$ and in the case $g =12$ it turns out to be an effective divisor in ${\mathcal R}_{12}$ of which we compute the cohomology class both  in ${\mathcal R}_{12}$ and in a partial compactification $\tilde{\mathcal R}_{12}$ introduced in \cite{fl}.

Recall that in \cite{fl} and in \cite{cefs}, using the existence of certain effective divisors on the moduli spaces ${\mathcal R}_{g,l}$ parametrizing isomorphism classes of pairs $[C,A]$ where $C$ is a smooth projective curve of genus $g$ and $A \in Pic^0(A)$ is a torsion line bundle of order $l$ (so ${\mathcal R}_g={\mathcal R}_{g,2}$), they prove that ${\mathcal R}_g= {\mathcal R}_{g,2}$ is a variety of general type if $g> 13$, and ${\mathcal R}_{g,3}$ is of general type if $g \geq 12$. Some of these divisors come from the study of the Koszul cohomology of the curve embedded by the linear systems $|\omega_C \otimes A|$, which has also been studied using  degeneration to binary curves in \cite{cf3} (see Theorem 0.6 of \cite{cefs} for the explicit computation of the classes). 
In section 6 we compare the class of our  divisor coming from the degeneracy locus of the Gaussian map $\mu_A$ with the divisor coming from the non vanishing of the Koszul cohomology group $K_{3,2}(\omega_C \otimes A)$ studied in \cite{fl} and \cite{cefs}. 


Finally we observe that both Gaussian maps of the form $\gamma^1_{\omega_C\otimes A, \omega_C \otimes A}$ (and also  of the form $\gamma^1_{\omega_C, \omega_C \otimes A}$ as in \cite{cv}), for a torsion line bundle $A \in Pic^0(A)$ of order $l> 2$  and Koszul cohomology can be studied by degeneration to binary curves with a small variation of the techniques used in this paper and in \cite{cf3}. One has to explicitly describe the line bundles of torsion $l >2$ in a similar way as we constructed the 2 torsion line bundles in section 2.

The paper is organized as follows: in section 2  we describe the construction given in \cite{cf} (and in \cite{cf3}) of Prym-canonical binary curves. In section 3 we recall the definition of Gaussian maps and we describe the first Gaussian map $\mu_A$ for Prym-canonical line bundles both for smooth curves and for Prym-canonical binary curves. In section 4 we prove by induction on the genus the surjectivity of the map $\mu_A$. In section 5 we describe the Maple script that gives the surjectivity of $\mu_A$ for $g =12$ and the injectivity statement for lower genus. In section 6 we compute the cohomology class in  $\tilde{\mathcal R}_{12}$ of the degeneracy locus of $\mu_A$. 

\section{Prym-canonical binary curves}

We first recall an explicit construction of binary curves embedded in ${\proj}^{g-1}$ by a Prym-canonical linear system $|\omega_C \otimes A|$ where $A^2 \cong \OO_C$, $A$ non trivial, which was given in \cite{cf}. Let $C$ be a binary curve of genus $g$, and $A\in Pic^0(C)$ a nontrivial line bundle. Then $H^0(C,\omega_C \otimes A)$ has dimension $g-1$
and  the restriction of $\omega_C \otimes A$ to  the component $C_j$ is
$\omega_{C_j}(D_j)$ where $D_j$ is the divisor of nodes on $C_j$. Since
$\omega_{C_j}(D_j)\cong \OO_{{\proj}^{1}}(g-1)$, the
components are embedded by a linear subsystem of
$\OO_{{\proj}^{1}}(g-1)$, hence they are projections from a point
of rational normal curves in ${\proj}^{g-1}$. Viceversa, let us
take two rational curves embedded in  ${\proj}^{g-2}$ by non
complete linear systems of degree $g-1$ intersecting transversally
at $g+1$ points. Then their union $C$ is a binary curve of genus
$g$ embedded either by a linear subsystem of $\omega_C$ or by a
complete linear system $|\omega_C \otimes A|$, where $A\in Pic^0(C)$
is nontrivial (see e.g. \cite{capo}, Lemma 10). Let us now recall the construction given in \cite{cf} (Lemma 3.1) of a
binary curve $C$ embedded in ${\proj}^{g-2}$ by  a linear system
$|\omega_C \otimes A|$ with $A^{\otimes 2}\cong \OO_C$, and $A$  non
trivial, and let us denote a binary curve with this embedding a Prym-canonical binary curve.

Assume that the first $g-1$ nodes, are $P_i=(0,...,0,1,0,...0)$ with 1 at the $i$-th place, $i=1,...,g-1$, the remaining two nodes are $P_g:=[t_1,...,t_{g-1}]$ with $t_i=0$ for $i
=1,...,[\frac{g}{2}]$, $t_i =1$, for $i =
[\frac{g}{2}]+1,...,g-1$.
and $P_{g+1}:=[s_1,...,s_{g-1}]$ with $s_i=1$ for $i
=1,...,[\frac{g}{2}]$, $s_i =0$, for $i =
[\frac{g}{2}]+1,...,g-1$.

Then the component $C_j$ is the image of the map
\begin{equation}\label{pcan}
\begin{gathered}\phi_j:{\proj}^1 \rightarrow {\proj}^{g-2}, \
j=1,2, \ \text{where} \\
\phi_1(t,u):= \left[\frac{tM_1(t,u)}{(t-a_{1,1}u)},..., \frac{tM_1(t,u)}{(t-a_{k,1}u)}, \frac{-M_1(t,u)d_1 a_{k+1,1}u}{A_1(t-a_{k+1,1}u)},..., \frac{-M_1(t,u)d_1a_{g-1,1}u}{A_1(t-a_{g-1,1}u)}\right]\\
\phi_2(t,u):= \left[\frac{tM_2(t,u)}{(t-a_{1,2}u)},..., \frac{tM_2(t,u)}{(t-a_{k,2}u)}, \frac{-M_2(t,u)d_2 a_{k+1,2}u}{A_2(t-a_{k+1,2}u)},..., \frac{-M_2(t,u)d_2a_{g-1,2}u}{A_2(t-a_{g-1,2}u)}\right]\\
\end{gathered}
\end{equation}
with $k :=[\frac{g}{2}]$, $M_j(t,u):= \prod_{r=1}^{g-1} (t-a_{r,j}u)$,  and
$A_j= \prod_{i=1}^{g-1} a_{i,j}$, $j=1,2$, $d_2$ is a nonzero constant and $d_1 = \frac{-d_2 A_1}{A_2}$.
Notice that we have $\phi_j([a_{l,j},1]) = P_l$, $l=1,...,g-1$, $\phi_j([0,1]) = P_g$, $\phi_j([1,0]) = P_{g+1}$, $j=1,2$.
In Lemma 3.1 of \cite{cf} it is proven  that for a general choice of $a_{i,j}$'s, $C=C_1\cup
C_2$ is a binary curve embedded in ${\proj}^{g-2}$ by  a linear
system $|\omega_C \otimes A|$ with $A^{\otimes 2}\cong \OO_C$ and $A$ nontrivial. In fact, recall
that $Pic^0(C) \cong {{\Co}^*}^g \cong {{\Co}^*}^{g+1}/{\Co}^*$, where ${\Co}^*$ acts diagonally,
and in Lemma 3.1 of \cite{cf} it is shown that our line bundle $A$ corresponds to the element $[(h_1,...,h_{g+1})] \in {{\Co}^*}^{g+1}/{\Co}^*$,
where $h_i=1$, for $i< [\frac{g}{2}]+1$, $h_i = -1$, for $i
=[\frac{g}{2}]+1,...,g-1$, $h_g=-1$, $h_{g+1} = 1$, so in particular $A$ is of
2-torsion.

\section{The first Prym-canonical Gaussian map}

\subsection{Gaussian maps}\label{sec:curva}

Let $Y$ be a smooth complex projective variety and let $\Delta_Y\subset
Y\times Y$ be the diagonal. Let $L$ and $M$ be line bundles on $Y$.
For a non-negative integer $k$, the \emph{k-th Gaussian map}
associated to these data is the restriction to diagonal map
\begin{equation}\label{gaussian1}\gamma^k_{L,M}:H^0(Y\times Y,I^k_{\Delta_Y}\otimes
L\boxtimes M )\rightarrow
H^0(Y,{I^k_{\Delta_Y}}_{|\Delta_Y}\otimes L\otimes M)\cong
H^0(Y,S^k\Omega_Y^1\otimes L\otimes M).
\end{equation}
The exact sequence
\begin{equation}
\label{Ik} 0 \rightarrow I^{k+1}_{\Delta_Y} \rightarrow
I^k_{\Delta_Y} \rightarrow S^k\Omega^1_Y \rightarrow 0,
\end{equation}(where $S^k\Omega^1_Y$ is identified to its image via the diagonal map), twisted by $L\boxtimes M$, shows that the domain of the $k$-th
Gaussian map is the kernel of the previous one:
$$\gamma^k_{L,M}:
ker \gamma^{k-1}_{L,M}\rightarrow H^0(S^k\Omega_Y^1\otimes L\otimes
M).$$

In this paper, we will exclusively deal with the first Gaussian
map for curves $C$, assuming also that $L =M$. 

The map $\gamma^0_L$ is the multiplication map of global
sections
$$H^0(C,L)\otimes
H^0(C,L)\rightarrow H^0(C,L^{\otimes 2})$$
 which obviously
vanishes identically on $\wedge^2 H^0(L)$.
 Consequently, $H^0(C
\times C, I_{\Delta_C}\otimes L\boxtimes L)$ decomposes as
$\wedge^2 H^0(L)\oplus I_2(L)$, where $I_2(L)$ is the kernel of
$S^2H^0(C,L)\rightarrow H^0(C,L^{\otimes 2})$. Since $\gamma^1_L$ vanishes
on symmetric tensors, one writes
$$\gamma^1_L:\wedge^2H^0(L)\rightarrow H^0(\Omega^1_C\otimes
L^{\otimes 2}).$$ 
Assume now that the line bundle $L$ is $\omega_C \otimes A$, with $A \in Pic^0(C)[2]$, and denote by
 \begin{equation}
 \label{mua}
 \mu_A:= \gamma^1_{\omega_C \otimes A}: \Lambda^2 H^0(  \omega_C \otimes A) \to H^0( \omega_C^{\otimes 3})
 \end{equation}
 the first Gaussian map. 
 
 It is useful to recall also the local definition of $\mu_A$.
 Given two sections $\sigma_i \in H^0(  \omega_C \otimes A)$, $i=1,2$, assume that locally $\sigma_i = f_i(z) dz \otimes l$, where $l$ is a local generator of $A$. The local expression of the first Gaussian map 
 $\mu_A: \Lambda^2 H^0(  \omega_C \otimes A) \rightarrow H^0( \omega_C^{\otimes 3}),$ is $\mu_A(\sigma_1 \wedge \sigma_2) := (f_1 f'_2 - f_2 f'_1) (dz)^3.$

Denote by  ${\mathcal R}_g^0$ the open subset of ${\mathcal R}_g$ where there exists the universal family $f:{\mathcal X} \rightarrow   {\mathcal R}_g^0$. If $b \in {\mathcal R}^0_g$, we have
$f^{-1}(b) = [C_b, A_b]$ where $C_b$ is a smooth irreducible curve of genus $g$ and $A_b \in Pic^0(C_b)[2]$ is a
line bundle of order $2$ on $C_b$. Denote by ${\mathcal P} \in Pic({\mathcal X})$ the corresponding Prym
bundle  and consider the Hodge bundle $f_*\omega_{f}$, where $\omega_f$ is the relative dualizing sheaf of $f$.

The maps $\mu_A$ glue together to give a map of vector bundles on ${\mathcal R}_g^0$,
\begin{equation}
\label{mu}
\mu: \Lambda^2(f_*(\omega_f \otimes {\mathcal P})) \to f_*((\omega_f \otimes {\mathcal P})^{\otimes 2} \otimes \omega_f) \cong f_*(\omega_f^{\otimes 3})
\end{equation}

\subsection{Strategy of the proof of surjectivity}

In the next sections we will prove the surjectivity
of the first Prym-canonical Gaussian map $\mu_A$ for the general point $[C,A] \in
{\mathcal R}_g$. We will do it by degeneration to binary curves
following the method used in \cite{ccm} for the second Gaussian
map of the canonical line bundle and also in \cite{cf} for the second Gaussian
map of Prym-canonical line bundles. Recall that ${\mathcal R}_g$
admits a suitable compactification $\overline{{\mathcal R}}_g$,
which is isomorphic to the coarse moduli space of the stack
${\bf{R}}_g$ of Beauville admissible double covers (\cite{b},
\cite{acv}) and to the coarse moduli space of the stack of Prym
curves (\cite{bcf}).

Consider the partial compactification $\tilde{\mathcal R}_g$ of ${\mathcal R}_g$ introduced in \cite{fl}. It is defined as the preimage through the forgetful map $\overline{\mathcal R}_g \to \overline{\mathcal M}_g$ of the open subvariety of $\overline{\mathcal M}_g$ consisting of 1-nodal irreducible stable curves of genus $g$. Denote by $\psi: {\mathcal X} \rightarrow \tilde{{\bf{R}}}_g$ the universal family  and by ${\mathcal P} \in Pic({\mathcal X})$ the corresponding Prym bundle as in \cite{fl} 1.1.
The map of vector bundles over ${\mathcal R}^0_g$, $\mu: \Lambda^2(f_*(\omega_f \otimes {\mathcal P})) \to f_*((\omega_f \otimes {\mathcal P})^{\otimes 2} \otimes \omega_f) \cong f_*(\omega_f^{\otimes 3})$ defined in \eqref{mu}, extends to a map
\begin{equation}
\label{mutilde}
\tilde{\mu}:  \Lambda^2(\psi_*(\omega_\psi \otimes {\mathcal P})) \rightarrow \psi_*( (\omega_\psi\otimes {\mathcal P})^{\otimes 2}\otimes {\Omega^1_\psi})\cong \psi_*(\omega_\psi^{\otimes 3}\otimes {\mathcal P}^{\otimes 2}\otimes {\mathcal I}_Z),
\end{equation}
where $Z=Sing(\psi)$, $\Omega^1_\psi \cong \omega_\psi
\otimes {\mathcal I}_Z$.

If $[C,A]$ is a point in $\tilde{\mathcal R}_g$, the local expression of
\begin{equation}
\label{gauss}
\mu_A: \Lambda^2H^0(\omega_C \otimes A) \rightarrow H^0((\omega_C \otimes A)^{\otimes 2} \otimes {\Omega^1_C})
\end{equation}
 is as follows.  In local coordinates let $\sigma_i = f_i(z) \xi \otimes l$, $i=1,2$ be two sections of $ \omega_C \otimes A$ where $\xi$ and $l$ are local generators of the line bundles $\omega_C$, respectively $A$, and define 
 \begin{equation}
 \label{locnodal}
 \mu_A(\sigma_1 \wedge \sigma_2) := (f_1 df_2 - f_2 df_1) \xi^2 \otimes l^2.
 \end{equation}
If $[C,A]$ is a Prym-canonical binary curve, the map $\mu_A: \Lambda^2H^0(\omega_C \otimes A) \rightarrow H^0(\omega_C^{\otimes 2} \otimes {\Omega^1_C})$ can be defined as in  \eqref{locnodal} and we will give an explicit description of it for the Prym-canonical binary  curves constructed in section 2.1. 

To prove by semicontinuity the surjectivity of $\mu_A$ for the general point in ${\mathcal R}_g$ in the following we will exhibit  Prym-canonical binary curves $[C,A]$  for which $\mu_A$ is surjective.

\subsection{1st Gaussian map for Prym-canonical binary curves}
Assume now that $[C,A]$ is a stable Prym-canonical binary curve as in section 2.1. The map $\mu_A: \Lambda^2H^0(\omega_C \otimes A) \rightarrow H^0(\omega_C^{\otimes 2} \otimes {\Omega^1_C})$ is defined as in \eqref{locnodal} and in the following we will give an explicit description of it.
To this end, let us first describe the space $H^0(\Omega^1_C \otimes \omega_C^{\otimes 2})$ for a Prym-canonical binary curve $C = C_1 \cup C_2$. If we denote by $\pi: \tilde{C} \to C$ the normalization map,  we have an exact sequence 
$$
0 \to {\mathcal T} \to \Omega^1_C \to {\mathcal F}_C \to 0, 
$$
where  ${\mathcal F}_C \cong \bigoplus_{i=1,2} \pi_* \omega_{C_i}$, ${\mathcal T}$ is a torsion sheaf supported at the nodes of $C$ and if $p$ is a node and the local equation of $C$ around $p$ is $xy=0$, ${\mathcal T}_p$ is a one dimensional vector space spanned by $ydx$. So, tensoring with $\omega_C^{\otimes 2}$ and taking global sections, we get an exact sequence
\begin{equation}
\label{codominio}
0 \rightarrow T \rightarrow H^0(\Omega^1_C \otimes \omega_C^{\otimes 2}) \rightarrow \bigoplus_{j=1,2} H^0(C_j, \omega_{C_j}^{\otimes 3}(2D_j)) \rightarrow 0,
\end{equation}
where, as above, $D_j$ is the divisor of nodes on $C_j$ and $T= H^0(C, {\mathcal T} \otimes \omega_C^{\otimes 2}) \cong \Co^{g+1}.$

Consider first the non torsion part of $\mu_A$, namely composition of $\mu_A$ with the projection on the space $\bigoplus_{j=1,2} H^0(C_j, \omega_{C_j}^{\otimes 3}(2D_j))$.  So, we can decompose the map $\nu_A= \nu_1 \oplus \nu_2$ where 
$$\nu_j: \Lambda^2 H^0(C, \omega_C \otimes A) \to H^0(C_j, \omega_{C_j}^{\otimes 3}(2D_j)) \cong H^0({\proj}^1, \OO_{{\proj}^1}(2g-4)).$$

In the coordinate $t$, the embeddings defining the two components of $C$ are
$\phi_j(t,1) = [\alpha_{1,j}(t)\ldots \alpha_{g-1,j}(t)]$, where the components $\alpha_{i,j}$ are defined as in \eqref{pcan}.
Let $\sigma_1 \ldots \sigma_{g-1}$ be the  basis of $H^0(\LB)$ given by the coordinate hyperplane sections. We have 
\begin{equation}
 \label{eq:polinu1} 
\nu_h (\sigma_{i} \wedge \sigma_{j}) = (\alpha_{i,h} \alpha_{j,h}' - \alpha_{j,h} \alpha_{i,h}')(\ud t)^3, \ h =1,2
\end{equation}

So, as an element of  $H^0({\proj}^1, \OO_{{\proj}^1}(2g-4))$ we can identify $\nu_h (\sigma_{i} \wedge \sigma_{j})$ with the polynomial $\nu_{ij,h}$, for which we have the following expressions, if in \eqref{pcan} we choose $d_2 =1$ and so $d_1 = -\frac{A_1}{A_2}$:

\begin{minipage}[t]{0.05\textwidth}
$ $ \\
\begin{equation}
\label{eq:nu}
\end{equation}
\end{minipage}
\begin{minipage}[t]{0.9\textwidth} 
\begin{align*}
\nu_{ij,h}(t)= &(a_{i,h} - a_{j,h}) t^2 \frac{ M_h^2(t,1) }{(t-a_{i,h})^2(t-a_{j,h})^2} &&  \quad\mbox{for } i<j \leq \left[\frac{g}{2}\right] \\ 
\nu_{ij,h}(t)=  & \frac{(a_{i,h}-a_{j,h})a_{i,h} a_{j,h}}{A_2^2}\frac{M_h^2(t,1)}{(t-a_{i,h})^2(t-a_{j,h})^2} && \quad\mbox{for } \left[\frac{g}{2}\right]+1 \leq i < j \\
\nu_{ij,h}(t)=  &\frac{(-1)^h a_{j,h}}{ A_2} (t^2- 2 a_{i,h} t + a_{j,h} a_{i,h})  \frac{M_h^2(t,1)}{(t-a_{i,h})^2(t-a_{j,h})^2} && \quad\mbox{for } i \leq \left[\frac{g}{2}\right] < j
\end{align*} 
\end{minipage}

We have the following commutative diagram with exact rows

\begin{equation}
\xymatrix{
0 \ar[r] & T \ar[r]^-{i}    & H^0(\varOmega_C^1 \otimes \omega_C^{\otimes 2}) \ar[r]^-{p=p_1 \oplus p_2} \ar@/_1pc/[l]_-{s}  & \bigoplus_{j=1,2} H^0({C}_j, \omega_{{C}_j}^{\otimes 3} (2 D_j) ) \ar[r] & 0 \\
0 \ar[r] & Ker\,(\nu_A) \ar[r] \ar[u]^{t} & {\bigwedge\nolimits}^2 H^0(C, \omega_C \otimes A)\ar[u]^{\mu_A} \ar[ur]_-{\nu_A=\nu_1 \oplus \nu_2}&
}
\end{equation}

where, if the local expression of $C$ in a neighborhood of the node $P_h$ is $xy=0$ and $\sigma \in  H^0(\varOmega_C^1 \otimes \omega_C^{\otimes 2})$, the $h$-th component of  $s(\sigma)$ at  the node $P_h$, $h =1,...,g$ is given by the coefficient of the torsion term $ydx$ in the local expression of $s(\sigma)$ at the node $P_h$. 
The map $t$ is the restriction to the kernel of $\nu_A$ of the map $\tau = s \circ \mu_A: \Lambda^2 H^0(C, \omega_C \otimes A) \to T, $ whose component at the node $P_h$, $h = 1,...,g$ can be easily computed  as follows:
\begin{equation}\label{eq:torsion}
\tau(\sigma_i \wedge \sigma_j)_h = \alpha_{j,1}'(a_{h,1}) \alpha_{i,2}'(a_{h,2}) - \alpha_{i,1}'(a_{h,1}) \alpha_{j,2}'(a_{h,2}).
\end{equation}
Similarly, for the component at the node $P_{g+1}$, one has:
\begin{equation}
\label{eq:torsionbis}
 \tau(\sigma_i \wedge \sigma_j)_{g+1} = g_{j,1}'(0) g_{i,2}'(0) - g_{i,1}'(0) g_{j,2}'(0),
 \end{equation}
where $g_{i,r}(u)  $ is the $i$-th component of $\phi_{r}(1,u)$, $r =1,2$ (see \eqref{pcan} for the explicit expression).
So we have $\mu_A = \nu_A \oplus \tau$, and in the next section we will prove by induction on the genus that it is surjective for the general Prym-canonical binary curve.

\section{Surjectivity}

Let $C \subset {\proj}^{g-2}$ be a Prym-canonical binary curve embedded by $\omega_C \otimes A$, with $A^{\otimes 2} \cong \OO_C$, as in \ref{sec:curva} and set $k :=[\frac{g}{2}]$ and denote by $\tilde{C}_r$ the partial normalization of $C$ at the node $P_r$ with $r= k$ if $g =2k$,  $r =k+1$ if $g = 2k+1$ and let $p_1,p_2$ the two points of $\tilde{C}_r$ over $P_r$. We make this choice of the node in order to obtain the Prym-canonical model for the curve $\tilde{C}_r$. In fact, observe that in this way, for a general choice of the $a_{i,j}$'s, the
projection  from $P_r$ sends the curve $C$ to the Prym-canonical
model of $\tilde{C}_r$ in ${\proj}^{g-3}$ given by the line bundle
$\omega_{\tilde{C}_r} \otimes  A'_r$  where $A'_r$
corresponds to the point $(h'_1,...,h'_{g-1}, 1) \in {{\Co}^*}^{g}/{\Co}^*$, with $h'_i = 1$
for $i\leq [\frac{g-1}{2}]$, $h'_i = -1$ for
$i=[\frac{g-1}{2}]+1,...,g-1$,  as described in section 2. In
fact  $[\tilde{C}_r,A'_r]$ is parametrized by
$a'_{i,j} = a_{i,j}$ for $i
\leq r-1$, $j =1,2$, $a'_{i,j} = a_{i+1,j}$ for $i
\geq r$, $j =1,2.$ So if we set $d'_{j} := \frac{d_j}{a_{r,j}}$, $j =1,2$, we clearly have a pair $[\tilde{C}_r,A'_r]$ as in \eqref{pcan}. For simplicity let us choose as above $d_2 =1$, so $d_1 = -\frac{A_1}{A_2}$, hence $d'_{2} := \frac{1}{a_{r,2}}$, $d'_{1} := -\frac{A_1}{A_2 a_{r,1}}$.

For the inductive step, we will use the same strategy used in \cite{ccm} and \cite{cf}.
Looking at the torsion and non-torsion part of the map $\mu_A$ separately, we have the following commutative diagrams with horizontal exact sequences
\begin{equation}\label{chi} 
\xymatrix{
0 \ar[r]	 & \bigoplus_{i=1,2} H^0(C_i, \omega_{C_i}^{\otimes 3} (2 \tilde{D}_i) ) \ar[r] & \bigoplus_{i=1,2} H^0({C}_i, \omega_{{C}_i}^{\otimes 3} (2 D_i) ) \ar@{->>}[r] & \bigoplus_{i=1,2}\OO_{2p_i} \\
0 \ar[r] & \bigwedge\nolimits^2 H^0(\tilde{C}_r,\omega_{\tilde{C}_r} \otimes A'_r) \ar[u]^{\tilde{\nu}_A} \ar[r] & \bigwedge\nolimits^2 H^0(C,\omega_C \otimes A) \ar[u]^{\nu_A} \ar[ur]^{\chi}
} 
\end{equation}

\begin{equation}
\label{tau} \xymatrix{
&0\ar[r] & \tilde T \ar[r] & T  \ar@{->>}[r] & T_P \\
&0 \ar[r] &\bigwedge\nolimits^2 H^0(\tilde{C}_r,\omega_{\tilde{C}_r} \otimes A'_r) \ar[u]^{\tilde\tau}  \ar[r]&  \bigwedge\nolimits^2 H^0(C,\omega_C \otimes A) \ar[u]^{\tau}
\ar[u]^{\tau} \ar[ur]_{\tau_P} & }
\end{equation}

where $P = P_r$, $D_i$ is the divisor of nodes of $C$ on $C_i$ and $\tilde{D_i} = D_i - p_i$ and $\nu_A, \tau$ and $\tilde\nu_A, \tilde\tau$ are the maps defined in the previous section for $C$ and $\tilde C_r$.
Hence, if  the maps $\tilde\nu_A \oplus \tilde{\tau}$ and $\chi \oplus \tau_P$ are surjective, then $\nu_A \oplus \tau = \mu_A$
 is also surjective. \\
Recall that the map $\nu_A$ is $\nu_1 \oplus \nu_2$, where $\nu_1$ and $\nu_2$ are defined in \eqref{eq:polinu1}, so we can write $\chi = \chi_1 \oplus \chi_2$, where $\chi_h$ is the composition of $\nu_i$ with the restriction to $\OO_{2p_i}$, $i = 1, 2$. 
In local coordinates around $P$, $\chi_h(\sigma_i \wedge \sigma_j)$ is the pair $(\nu_{ij,h}(a_{r,h}), \nu_{ij,h}'(a_{r,h}))$ (where $r = k$ for $g=2k$ and $r = k + 1$ for $g=2k+1$) corresponding to the evaluation of the polynomial $\nu_{ij,h}(t)$ of \eqref{eq:polinu1} and of its derivative at $P$.
Instead $\tau_P$ is nothing else than the torsion part of the map $\mu_A$ in the node $P$ (as in \eqref{eq:torsion}).
So we can prove the following

\begin{TEO}\label{thm20}
If $[C=C_1\cup C_2, \omega_C\otimes A]$  is a Prym-canonical general binary curve of genus $g \geq 12$, then $\mu_A$ is surjective for $C$.
\end{TEO}
\proof
The case $g=12$ is done by computations with Maple (see section 5).
For $g>12$ we proceed by induction, as explained above, so we assume that the map $\tilde{\nu}_A \oplus \tilde{\tau}$ of the above diagrams is surjective, hence we only have to check the maximality of the rank of the following $5 \times \frac{1}{2}(g-1)(g-2)$ matrix, for every $g \geq 13$:
\begin{equation*}
\begin{pmatrix}
\nu_{ij,1} |_{t=a_{k,1}} \\
\frac{\ud}{\ud t} (\nu_{ij,1}) |_{t=a_{k,1}} \\
\nu_{ij,2} |_{t=a_{k,2}} \\
\frac{\ud}{\ud t} (\nu_{ij,2}) |_{t=a_{k,2}} ) \\
\tau_{ij}(P_k)
\end{pmatrix}
\quad \boxed{\text{if } \mathbf{g=2k}} \qquad
\begin{pmatrix}
\nu_{ij,1} |_{t=a_{{k+1},1}} \\
\frac{\ud}{\ud t} (\nu_{ij,1}) |_{t=a_{{k+1},1}} \\
\nu_{ij,2} |_{t=a_{{k+1},2}} \\
\frac{\ud}{\ud t} (\nu_{ij,2}) |_{t=a_{{k+1},2}} ) \\
\tau_{ij}(P_{k+1}) 
\end{pmatrix}
 \quad \boxed{\text{if } \mathbf{g=2k+1}}
\end{equation*}
where $1 \leq i < j \leq g-1$ are the columns indexes. \\

Let us start with the even case $g=2k$.
We will see that the $5 \times 5$ submatrix with column indexes $(i,j)=(1,k),(2,k),(k,g-2),(k,g-1),(k-1,k+1)$ is invertible for some appropriate choices of the parameters $a_{\ell,\epsilon}$, $\epsilon=1,2$. \\
For $\boxed{i=1,2; \: j=g-2,g-1; \: \epsilon=1,2}$, taking $\boxed{\delta_{\epsilon}=(-1)^{\epsilon+1}}$, we have (see equations \eqref{eq:nu})
\begin{align*}
&N_{i,\epsilon} := \nu_{ik,\epsilon}|_{t=a_{k,\epsilon}} = a_{k,\epsilon}^2 (a_{i,\epsilon}-a_{k,\epsilon}) \prod_{\ell \neq i,k}(a_{k,\epsilon}-a_{\ell,\epsilon})^2 \\
&N_{i,\epsilon}' := \frac{\ud (\nu_{ik,\epsilon})}{\ud t}|_{t=a_{k,\epsilon}} = 2 a_{k,\epsilon} (a_{i,\epsilon}-a_{k,\epsilon}) \prod_{ \ell \neq i,k}(a_{k,\epsilon}-a_{\ell,\epsilon}) \left( \prod_{ \ell \neq i,k}(a_{k,\epsilon}-a_{\ell,\epsilon}) + a_{k,\epsilon} \sum_{h \neq i,k} \prod_{ \ell \neq i,k,h}(a_{k,\epsilon}-a_{\ell,\epsilon}) \right) \\
&N_{j,\epsilon} := \nu_{kj,\epsilon}|_{t=a_{k,\epsilon}} = \delta_{\epsilon} \frac{a_{k,\epsilon} a_{j,\epsilon}}{A_2} (a_{k,\epsilon}-a_{j,\epsilon}) \prod_{\ell \neq k,j}(a_{k,\epsilon}-a_{\ell,\epsilon})^2 \\
&N_{j,\epsilon}' := \frac{\ud (\nu_{kj,\epsilon})}{\ud t}|_{t=a_{k,\epsilon}}  =\delta_{\epsilon} \frac{2 a_{k,\epsilon} a_{j,\epsilon}}{A_2} (a_{k,\epsilon}-a_{j,\epsilon}) \prod_{\ell \neq k,j}(a_{k,\epsilon}-a_{\ell,\epsilon}) \cdot \sum_{h \neq k,j} \prod_{\ell \neq k,j,h}(a_{k,\epsilon}-a_{\ell,\epsilon}),
\end{align*}
while if $i,j \neq k$, we have $  \nu_{ij,\epsilon}|_{t=a_{k,\epsilon}} =0$ and $\frac{\ud (\nu_{ij,\epsilon})}{\ud t}|_{t=a_{k,\epsilon}} =0$, hence every entry of our matrix with $i,j \neq k$ is zero.
So, it's sufficient to show that, for an appropriate  choice of the parameters $a_{i,\epsilon}$, we have   $\tau_{k-1,k+1}(P_k) \neq 0$ and  the following $4\times4$ matrix is invertible:
\begin{equation}\label{matrice4}
\begin{pmatrix}
N_{1,1}  & N_{2,1}  & N_{g-2,1}  & N_{g-1,1} \\
N_{1,1}'  & N_{2,1}'  & N_{g-2,1}'  & N_{g-1,1}' \\
N_{1,2}  & N_{2,2}  & N_{g-2,2}  & N_{g-1,2} \\
N_{1,2}'  & N_{2,2}'  & N_{g-2,2}'  & N_{g-1,2}'
\end{pmatrix}
\end{equation}

First of all, perform the following simplifications:
\begin{itemize}
\item Multiply the last two columns by $A_2$.
\item Divide the first row by $ a_{k,1} \prod_{\substack{\ell \neq 1,2,k \\ \ell \neq g-2,g-1}}(a_{k,1}-a_{\ell,1})^2 \cdot (a_{k,1}-a_{1,1})(a_{k,1}-a_{2,1})(a_{k,1}-a_{g-2,1})(a_{k,1}-a_{g-1,1})$ and similarly for the third row.
\item Divide the second row by $2 a_{k,1} \prod_{\ell \neq k }(a_{k,1}-a_{\ell,1}) $  and similarly for the fourth row.
\end{itemize}
We obtain
\begin{align*}
&N_{i,\epsilon} = - a_{k,\epsilon} \prod_{\substack{\ell = 1,2,g-2,g-1 \\ \ell \neq i}}(a_{k,\epsilon}-a_{\ell,\epsilon}) \\
&N_{i,\epsilon}' = - \left( \prod_{ \ell \neq i,k}(a_{k,\epsilon}-a_{\ell,\epsilon}) + a_{k,\epsilon} \sum_{h \neq i,k} \prod_{ \ell \neq i,k,h}(a_{k,\epsilon}-a_{\ell,\epsilon}) \right)  \\
&N_{j,\epsilon} = \delta_{\epsilon} a_{j,\epsilon} \prod_{\substack{\ell = 1,2,g-2,g-1 \\ \ell \neq j}}(a_{k,\epsilon}-a_{\ell,\epsilon}) \\
&N_{j,\epsilon}' = \delta_{\epsilon} a_{j,\epsilon} \sum_{h \neq k,j} \prod_{\ell \neq k,j,h}(a_{k,\epsilon}-a_{\ell,\epsilon})
\end{align*}

Now, let $p_\epsilon=\prod_{\substack{\ell \neq 1,2,k \\ \ell \neq g-2,g-1}}(a_{k,\epsilon}-a_{\ell,\epsilon})$, $p_{h,\epsilon}=\frac{p_\epsilon}{(a_{k,\epsilon}-a_{h,\epsilon})} $ and $\hat{a}_{\ell,\epsilon}=(a_{k,\epsilon}-a_{\ell,\epsilon})$, for all $\ell \neq k$, so that we can write
\[ N_{i,\epsilon}' = - \left( p_\epsilon (\prod_{\substack{\ell = 1,2,g-2,g-1 \\ \ell \neq i}}\hat{a}_{\ell,\epsilon}) + a_{k,\epsilon} (\prod_{\substack{\ell = 1,2,g-2,g-1 \\ \ell \neq i}}\hat{a}_{\ell,\epsilon} )(\sum_{\substack{h \neq 1,2,k, \\ g-2,g-1}} p_{h,\epsilon} )+ a_{k,\epsilon} p_\epsilon \left( \sum_{\substack{n = 1,2,g-2,g-1 \\ n \neq i}}\frac{\prod_{\substack{\ell = 1,2,g-2,g-1 \\ \ell \neq i}}\hat{a}_{\ell,\epsilon}}{\hat{a}_{n,\epsilon}} \right) \right) \] 

\[N_{j,\epsilon}' = \delta_{\epsilon} \left( a_{j,\epsilon} (\prod_{\substack{\ell = 1,2,g-2,g-1 \\ \ell \neq j}}\hat{a}_{\ell,\epsilon})( \sum_{\substack{h \neq 1,2,k, \\ g-2,g-1}} p_{h,\epsilon}) + a_{j,\epsilon} p_\epsilon \left( \sum_{\substack{n = 1,2,g-2,g-1 \\ n \neq j}}\frac{\prod_{\substack{\ell = 1,2,g-2,g-1 \\ \ell \neq j}}\hat{a}_{\ell,\epsilon}}{\hat{a}_{n,\epsilon}} \right) \right) \] 
Now we choose the following values for the parameters:
\begin{equation}\label{sceltapar}
\boxed{a_{i,1}=i \cdot a, \quad a \neq 0,1 \qquad a_{i,2}=i}
\end{equation}

So, we obtain $\sum_{\substack{h \neq 1,2,k, \\ g-2,g-1}} p_{h,\epsilon} = 0$ for $\epsilon = 1,2$ and $p_1 = a^{g-6} p_2 = a^{g-6} (-1)^{k-3} ((k-3)!)^2$.
Dividing the first row by  $a^4$,  the second row by $a^3p_1$,  the fourth row by $p_2$, and the first and third row by $(k-1)(k-2)$, we obtain the matrix 
\[ \begin{pmatrix}
-k(k-2)        & -k(k-1)         & -2(k-1)^2    & -(2k-1)(k-2)  \\
(k-2)^2      & 2(k-1)^2     & -2(k-1)^3    & -(2k-1)(k-2)^2 \\
-k(k-2)         & -k(k-1)         & 2(k-1)^2   & (2k-1)(k-2)  \\
(k-2)^2      & 2(k-1)^2      & 2(k-1)^3   & (2k-1)(k-2)^2
\end{pmatrix} \] 
that is invertible.
Finally we can check that with the parameters \eqref{sceltapar} we have
\[  \tau_{k-1,k+1}(P_k) =  \frac{2k(k+1)a^{g-2}}{A_2} \prod_{\ell \neq k-1,k+1,k} (k-\ell)^2   \neq 0 \]

The computations and strategy for the odd case $g=2k+1$ are similar. We will prove that with the same choice for the $a_{i,\epsilon}$'s as in the even case, the  $5 \times 5$ matrix with column indexes $(2,k+1),(3,k+1),(k+1,g-2),(k+1,g-1),(k-1,k+2)$ is invertible. 

For the indexes $(i,j)=(2,k+1),(3,k+1),(k+1,g-2),(k+1,g-1)$ and $\delta_{\epsilon}=(-1)^{\epsilon+1}$, we have
$$
N_{i,\epsilon} = \nu_{i\,k+1,\epsilon}|_{t=a_{k+1,\epsilon}} = - \delta_{\epsilon} \prod_{\ell \neq i,k+1}(a_{k+1,\epsilon}-a_{\ell,\epsilon})^2    \cdot    \frac{ (a_{k+1,\epsilon}^2- a_{i,\epsilon} a_{k+1,\epsilon}) a_{k+1,\epsilon} }{A_2}$$
$$N_{i,\epsilon}' = \frac{\ud (\nu_{ik+1,\epsilon})}{\ud t}|_{t=a_{k+1,\epsilon}} =- \delta_{\epsilon} \frac{2a_{k+1,\epsilon}(a_{k+1,\epsilon}-a_{i,\epsilon}) \prod_{\ell \neq i,k+1}(a_{k+1,\epsilon}-a_{\ell,\epsilon})}{A_2}  \cdot  [\prod_{\ell \neq i,k+1}(a_{k+1,\epsilon}-a_{\ell,\epsilon}) +$$
$$+a_{k+1,\epsilon}\sum_{h \neq i,k+1} \prod_{\ell \neq i,k+1,h}(a_{k+1,\epsilon}-a_{\ell,\epsilon})] $$
$$N_{j,\epsilon}  = \nu_{k+1\,j,\epsilon}|_{t=a_{k+1,\epsilon}} = \prod_{\ell \neq j,k+1}(a_{k+1,\epsilon}-a_{\ell,\epsilon})^2     \cdot     \frac{(a_{k+1,\epsilon}-a_{j,\epsilon}) a_{k+1,\epsilon} a_{j,\epsilon} }{A_2^2} $$
$$N_{j,\epsilon}' = \frac{\ud (\nu_{k+1j,\epsilon})}{\ud t}|_{t=a_{k+1,\epsilon}}= 2 \cdot \prod_{\ell \neq k+1,j}(a_{k+1,\epsilon}-a_{\ell,\epsilon}) \cdot \sum_{h \neq k+1,j} \prod_{\ell \neq k+1,j,h}(a_{k+1,\epsilon}-a_{\ell,\epsilon}) \cdot \frac{(a_{k+1,\epsilon}-a_{j,\epsilon}) a_{k+1,\epsilon} a_{j,\epsilon} }{A_2^2}
$$

Making analogous simplifications as in the even case, and choosing as above the parameters $a_ {i,1} = ia$, $a \neq 0,1$, $a_ {i,2} = i$, we obtain the matrix
\[ \begin{pmatrix}
-(k+1)(k-1) & -(k+1)(k-2) & -(k-2) & -(k-1)      \\
-(k^2-2k-1) & -(k^2-4k-2) & -(2k-2) & -(2k-1)    \\
(k+1)(k-1) & (k+1)(k-2) & -(k-2) & -(k-1)      \\
(k^2-2k-1) & (k^2-4k-2) & -(2k-2) & -(2k-1)     
\end{pmatrix} \]
that is an invertible. Finally we observe that if $i,j \neq k+1$, we have $  \nu_{ij,\epsilon}|_{t=a_{k+1,\epsilon}} =0$ and $\frac{\ud (\nu_{ij,\epsilon})}{\ud t}|_{t=a_{k+1,\epsilon}} =0$, and 
\[ \tau_{k-1,k+2}(P_{k+1}) = \frac{-4(k+1)(k+2)a^{g-2}}{A_2} \prod_{\ell \neq k-1,k+2,k+1} (k+1-\ell)^2  \neq 0. \]

\endproof


\section{Maple scripts for computations and injectivity for low genus}

In this section we want to show the following
\begin{TEO}\label{thmgenbassi}
If $[C=C_1\cup C_2, \omega_C\otimes A]$ is a Prym-canonical general binary curve of genus $3 \leq g \leq 12$, then $\mu_A$ is injective for $C$.
\end{TEO}

\proof

For the cases  $g=3,4$, assume that $C$ is a smooth curve and let $\omega_C \otimes A$ be a  Prym-canonical line bundle, then, given a non zero decomposable vector  $\sigma_1 \wedge \sigma_2 \in \Lambda^2H^0(\LB)$, we have $\mu_A(\sigma_1 \wedge \sigma_2)\neq 0$, as one can easily check, or find in  \cite{wah92}, (1.4). This concludes the proof for genus 3 curves. 
Now, by a simple linear algebra remark (see e.g. \cite{cf1}, (2.4)), one sees that $rank(\mu_A) \geq dim(G(2,h^0(\omega_C \otimes A) ))+ 1 = 2 h^0(\omega_C \otimes A) -3$, where $G(2,h^0(\omega_C \otimes A) )$ is the Grassmannian of 2-dimensional subspaces of $H^0(\omega_C \otimes A)$. 
So for any curve of genus 4, we have $rank(\mu_A) \geq 2 h^0(\omega_C \otimes A) -3= dim( \Lambda^2 H^0(\omega_C \otimes A))=3$.

For $g \geq 5$ we can use the model of Prym-canonical binary curves presented in section 2.
We will choose some random values for the parameters $a_{i,j}$ and check, making the computation with Maple, that there exists  a Prym-canonical curve constructed like in section 2 for which the map $\mu_A$ has maximal rank.
Observe that for $g=12$ we have $dim(\Lambda^2 H^0(\LB))=dim(H^0(\varOmega^1_C \otimes \omega_C^{\otimes 2})=55$, so the maximality of the rank of $\mu_A$ for genus $12$ is the basis of the induction for the proof of  the surjectivity for $g \geq 12$, as explained in the previous section.
For $4 < g < 12$ we obtain, by semicontinuity, the injectivity of $\mu_A$ for a general point of $\mathcal{R}_g$.
The following is the Maple script used to construct the matrices of the torsion and non torsion part of $\mu_A$ in the way explained in section 3.3 and to check the rank for genus $4 \leq g \leq 12$:

\begin{verbatim}
(1)  with(LinearAlgebra):                                                       
     v:=[1,2,3,4,5,6,7,8,9,10,11]:  w:=[326,-28,-875,-97,20,-651,-523,-306,369,-31,99]:
     for g from 4 to 12 do
 	   unassign(a): unassign(b):  a[g]:=0: b[g]:=0:  k:=floor(g/2):
(5)    for i from 1 to k do P[i]:=0 end do:
       for i from k+1 to g-1 do P[i]:=1 end do:
       A2:=mul(b[i],i=1..g-1):  
       c[1]:=[seq(-P[i]*a[i]*A2,i=1..g-1)]:  c[2]:=[seq(P[i]*b[i]/A2,i=1..g-1)]:
       for i from 1 to k do delta[i]:=1 end do:
(10)   for i from k+1 to g-1 do delta[i]:=0 end do:
       M[1]:=mul(t-a[i], i = 1..g-1):  M[2] := mul(t-b[i], i = 1..g-1):
       for i from 1 to g-1 do 
         alpha[i,1]:=M[1]*(delta[i]*t-c[1][i])/(t-a[i]);  
         alpha[i,2]:=M[2]*(delta[i]*t-c[2][i])/(t-b[i])
       end do:                                                                  
(15)   for i from 1 to g-1 do  for h from 1 to 2 do  
         alphaprimo[i,h]:=diff(alpha[i,h],t);  
       end do:  end do:
       for i from 1 to g-2 do  for j from i+1 to g-1 do  for h from 1 to 2 do  
         R[h][i,j]:=alpha[i,h]*alphaprimo[j,h] - alpha[j,h]*alphaprimo[i,h];  
(20)   end do:  end do:  end do:                                                                  
       MM[1]:=mul(1-a[i]*u, i = 1..g-1):  MM[2] := mul(1-b[i]*u, i = 1..g-1):
       for i from 1 to g-1 do 
         tauGp1a[i]:=diff(MM[1]*(delta[i]-c[1][i]*u)/(1-a[i]*u),u); 
         tauGp1b[i]:=diff(MM[2]*(delta[i]-c[2][i]*u)/(1-b[i]*u),u)        
(25)   end do:
       for i from 1 to g-2 do  for j from i+1 to g-1 do  for h from 1 to g do
         tia:=eval(alphaprimo[i,1],t=a[h]);  tja:=eval(alphaprimo[j,1],t=a[h]);                                 
         tjb:=eval(alphaprimo[j,2],t=b[h]);  tib:=eval(alphaprimo[i,2],t=b[h]);
         tau[h][i,j] := (tja)*(tib) - (tia)*(tjb)
(30)   end do:  end do:  end do:
       for i from 1 to g-2 do   for j from i+1 to g-1 do
         seqnu:=seq( seq(coeff(R[h][i,j],t,n), n=0..2*g-4), h=1..2 );       
         tGp1ia:=eval(tauGp1a[i],u=0);  tGp1jb:=eval(tauGp1b[j],u=0);
         tGp1ja:=eval(tauGp1a[j],u=0);  tGp1ib:=eval(tauGp1b[i],u=0);
(35)     seqtau:=seq(tau[h][i,j], h=1..g);                                   
         riga[i,j]:=[ seqnu, seqtau, (tGp1ja)*(tGp1ib)-(tGp1ia)*(tGp1jb) ] 
       end do:  end do:
       ttt[g]:=[ seq(seq(riga[i,j], j=i+1..g-1), i=1..g-2) ]:
       T[g]:=convert(ttt[g], Matrix):                                         
(40)   mappamu[g]:=eval( eval(T[g], a=v[1..g-1]), b=w[1..g-1] ):
       rango[g]:=Rank(mappamu[g]):
     end do:
     for g from 4 to 12 do T[g]; rango[g]; end do;
\end{verbatim}

We will  briefly explain the code.
At the second line we take the vectors with the choice of the values that we will use for parameters.
Form line 4 to line 14 we calculate the functions $\alpha_{i,j}$ which represent the embeddings $\phi_j(t,1)$.
Note that the sequence memorized in the variable \texttt{a} is the sequence of  the parameters $a_{i,1}$, $i=1,\ldots,g-1$, while the parameters $a_{i,2}$ are memorized in the variable \texttt{b}.
From line 15 to 20, we calculate $\alpha_{i,j}'$ and the polynomials $\nu_{ij,h}$ (\textit{i.e.}~\texttt{R[h][i,j]}) as in equation~\eqref{eq:polinu1}.
In lines 21--25 we calculate the functions $g_{i,r}'$ (called \texttt{tauGp1a} and \texttt{tauGp1b} in the code), \textit{i.e.}~the derivative of the $(g+1)$-coordinate of $\phi_r(1,u)$, $r=1,2$, needed for the torsion part in the last node (see equation~\eqref{eq:torsionbis}).
In lines 24--30 we calculate the $h$-component, for $h=1,\ldots,g$, of the torsion for every $\sigma_i \wedge \sigma_j$ (as in  equation~\eqref{eq:torsion}).
Then we take the coefficients of the  polynomials $\nu_{ij,h}$ (line 32) and assemble the matrix $\binom{\nu_A}{\tau}$ (line 39, for convenience in the code, we built the transpose of this matrix).
Finally, we replace the values in the parameters (line 40) and we calculate the rank of the matrix, which turns out to be maximal for $g=4,\ldots,12$. 
\endproof

\section{The class}

In the previous sections we have proved by semicontinuity that the
first Gaussian map 
$$\mu_A: \Lambda^2H^0(\omega_C \otimes A) \rightarrow H^0(\omega_C^{\otimes 2} \otimes {\Omega^1_C})$$
has maximal rank for the general pair
$[C,A]$ in ${\mathcal R}_{12}$. 
Notice that for $g =12$,
$dim( \Lambda^2H^0(\omega_C \otimes A)) = dim(H^0(\omega_C^{\otimes 2} \otimes {\Omega^1_C})) =
55. $ Consider the locus ${\mathcal D} = \{[C,A] \in {\mathcal
R}_{12} \ | \ rk(\mu_A) <55 \}$. We have proved that ${\mathcal
D} \neq {\mathcal R}_{12}$, hence it is an
effective divisor in ${\mathcal R}_{12}.$ Let $\pi:
\overline{\mathcal R}_{g} \rightarrow \overline{{\mathcal M}}_g$
be the finite map which extends the forgetful map ${\mathcal R}_g
\rightarrow {\mathcal M}_g$ (see \cite{fl} Section 1). The partial
compactification $\tilde{\mathcal R}_g$ of ${\mathcal R}_g$
introduced in \cite{fl} Section 1 is the inverse image
$\pi^{-1}(\tilde{\mathcal M}_g)$, where $\tilde{\mathcal
M}_g:={\mathcal M}_g\cup \tilde{\Delta }_0$ and $\tilde{\Delta
}_0$ is the locus of one-nodal irreducible curves. Denote by $\psi:
{\mathcal X} \rightarrow \tilde{{\bf{R}}}_g$ the universal family
and by ${\mathcal P} \in Pic({\mathcal X})$ the corresponding Prym
bundle as in \cite{fl} 1.1. Assume $g =12$, then if
$\tilde{{\mathcal D}}$ is the closure of ${\mathcal D}$ in
${\tilde{\mathcal R}_{12}}$, $\tilde{{\mathcal D}}$ is the
degeneracy locus of the map
$$  \Lambda^2(\psi_*(\omega_\psi \otimes {\mathcal P})) \rightarrow \psi_*( (\omega_\psi\otimes {\mathcal P})^{\otimes 2}\otimes {\Omega^1_\psi})\cong \psi_*(\omega_\psi^{\otimes 3}\otimes {\mathcal P}^{\otimes 2}\otimes {\mathcal I}_Z)
$$
of
\eqref{mutilde}, where $Z=Sing(\psi)$, $\Omega^1_\psi \cong \omega_\psi \otimes {\mathcal I}_Z$. Denote by ${\mathcal F}_i :=
\psi_{*}(\omega_\psi^{\otimes i} \otimes {\mathcal P}^{\otimes i})$.
Using Grothendieck-Riemann-Roch and Proposition 1.6 of \cite{fl}
one computes as in Proposition 1.7 of \cite{fl}
$$
 c_1({\mathcal F}_i) = \frac{i(i-1)}{2} (12 \lambda - \delta'_0- \delta''_0 - 2 \delta_0^{ram}) +
 \lambda - \frac{i^2}{4} \delta_0^{ram},
$$
where $\lambda$ is the pullback of the Hodge class $\lambda \in {\overline{\mathcal M}}_g$ and $ \delta'_0$, $\delta''_0$, and $\delta_0^{ram}$ are the boundary classes defined in \cite{fl} section 1.
So we have
$
 c_1({\mathcal F}_1) = \lambda - \frac{\delta_0^{ram}}{4}$, $c_1(\Lambda^2(\psi_*(\omega_\psi \otimes {\mathcal P})) = 10 \lambda -\frac{5}{2} \delta_0^{ram}$.

Notice that  by Grothendieck-Riemann-Roch we have
 $$c_1( \psi_*( \omega_\psi^{\otimes 3} \otimes {\mathcal P}^{\otimes 2} \otimes {\mathcal I}_Z)) =
 \psi_*[(1 + c_1( \omega_\psi^{\otimes 3} \otimes {\mathcal P}^{\otimes 2})
 + \frac{1}{2} c_1^2( \omega_\psi^{\otimes 3} \otimes {\mathcal P}^{\otimes 2})- [Z])
 \cdot (1 - \frac{c_1(\omega_\psi)}{2} + \frac{ c_1(\omega_\psi)^2+ [Z]}{12}) ]_2$$
$$= 37 \lambda -4 (\delta'_0 + \delta''_0) -9 \delta_0^{ram}, $$
since $\psi_*(c_1(\omega_\psi) \cdot {\mathcal P}) =0$, $\psi_*(c_1({\mathcal P})^2) = - \delta_0^{ram}/2$,
by Proposition 1.6 of \cite{fl} and by Mumford's formula, $\psi_*(c_1(\omega_\psi)^2) = 12 \lambda - \psi_*([Z])$
and $\psi_*[Z] = \delta'_0 + \delta''_0+ 2 \delta_0^{ram}$ (\cite{fl}, 1.1). So, finally we have
$$ c_1(\tilde{{\mathcal D}}) = c_1( \psi_*( \omega_\psi^{\otimes 3} \otimes {\mathcal P}^{\otimes 2}
\otimes {\mathcal I}_Z))\cdot rk(    \Lambda^2(\psi_*(\omega_\psi \otimes {\mathcal P})  )- c_1(\Lambda^2(\psi_*(\omega_\psi \otimes {\mathcal P}))\cdot
rk  (\psi_*( \omega_\psi^{\otimes 3} \otimes {\mathcal P}^{\otimes 2}\otimes {\mathcal I}_Z))= $$
$$=55 (27 \lambda -4 (\delta'_0 + \delta''_0) -\frac{13}{2} \delta_0^{ram}),$$
and $c_1({\mathcal D}) = 1485\lambda$, hence ${\mathcal D}$ is an effective divisor in
${\mathcal R}_{12}$, $\tilde{{\mathcal D}}$ is an effective divisor in $\tilde{{\mathcal R}}_{12}$
and if we denote by $\overline{{\mathcal D}}$ the closure of ${\mathcal D}$ in $\overline{\mathcal R}_{12}$,
we have computed
\begin{equation}
\label{Dbar}
c_1(\overline{{\mathcal D}}) = 55(27 \lambda -4 (\delta'_0 + \delta''_0) -\frac{13}{2} \delta_0^{ram}- ...)
\end{equation}
In fact, since the partial compactification $\tilde{{\mathcal
R}}_{g} \subset \overline{\mathcal R}_{g}$ has the property that
$\pi^{-1}({\mathcal M}_g \cup \Delta_0)-\tilde{{\mathcal R}}_{g} $
has codimension $\geq 2$, the expression (\ref{Dbar}) computes the
coefficients of $\lambda$, $\delta'_0$, $\delta''_0$,
$\delta_0^{ram}$ in $c_1(\overline{{\mathcal D}})$.

Unfortunately the ratio between the coefficient of $\lambda$ and the coefficient of $\delta'_0+ \delta''_0$ is $>13/2$, hence it does not seem possible to deduce anything on the Kodaira dimension of ${{\mathcal R}}_{12}$ using this divisor (see e.g. Remark 3.5 of \cite{cefs}).

Finally note that in \cite{cefs}  and in \cite{cf3} it is proven that for the general point of ${{\mathcal R}}_{12}$ the Koszul cohomology $K_{3,2}(C, \omega_C \otimes A)$ vanishes and in \cite{fl} (see also \cite{cefs}, Theorem 0.6) they compute the class of the divisor $Z$ given by the points $[C,A]$ such that the Koszul cohomology $K_{3,2}(C, K_C \otimes A) \neq 0$  in $\tilde{{\mathcal
R}}_{12}$, which is 
$$[\tilde{Z}]= 56(\frac{13}{2} \lambda -(\delta'_0 + \delta''_0) -\frac{3}{2} \delta_0^{ram}), $$
where $\tilde{Z}$ denotes the closure of $Z$ in $\tilde{{\mathcal
R}}_{12}$. 
So, since we have the following formula for the canonical class of $\tilde{{\mathcal R}}_{12}$, $$K_{ \tilde{{\mathcal R}}_{12}} = 13 \lambda -2(\delta'_0 + \delta''_0) -3\delta_0^{ram},$$
if we denote by $\overline{Z}$ the closure of $Z$ in $ \overline{{\mathcal R}}_{12}$,  using Remark 3.5 of \cite{cefs}, one immediately sees that  $K_{ \overline{{\mathcal R}}_{12}}-\frac{1}{28}[\overline{Z}]$  is an effective classe in $ \overline{{\mathcal R}}_{12}$, so $K_{ \overline{{\mathcal R}}_{12}}$ is an effective class.

\end{document}